\documentclass[12pt,reqno]{amsart}

\usepackage{amsmath,amssymb,amsfonts,amsthm}
\newtheorem{theorem}{Theorem}[]

\usepackage{tikz}
\usetikzlibrary{matrix,arrows.meta,positioning, calc}
\usepackage[a4paper,
top=24mm,bottom=24mm,left=26mm,right=26mm
]{geometry}
\usepackage{graphicx}
\usepackage{booktabs}
\usepackage{mathtools}
\usepackage{nicematrix}
\usepackage{makecell}

\usepackage[
]{cochineal}
\usepackage[cochineal,vvarbb,
cochf,cochrho]{newtxmath}
\usepackage[cal=zapfc, calscaled=1.15,
scr=boondox
]{mathalpha}

\usepackage[numbers,sort]{natbib}
\usepackage[
bookmarks=true,
bookmarksdepth=subsection,
bookmarksnumbered=true,
breaklinks=true,
colorlinks=true,
linkcolor=blue,
citecolor=magenta
]{hyperref}

\usepackage{listings}
\usepackage{xcolor}
\usepackage{tcolorbox}

\DeclareMathOperator{\rank}{rank}

\DeclareMathOperator{\Aut}{Aut}
\DeclareMathOperator{\Tr}{Tr}

\DeclareMathOperator{\Hom}{Hom}

\DeclareMathDelimiter{\Norm}{\mathord}{largesymbols}{"3E}{largesymbols}{"3E}

\allowdisplaybreaks
\begin{document}

\title[]{Character Variety of the Punctured Sphere
  and
  the Determinantal Variety}
\author[]{Kazuhiro Hikami}
\address{Faculty of Mathematics,
  Kyushu University,
  Fukuoka 819-0395, Japan.}

\email{
  \texttt{khikami@gmail.com}
}
\date{November 15, 2026}

\begin{abstract}
  We show that the $\mathrm{SL}(2;\mathbb{C})$ character variety of
  the $n$-punctured sphere is described via the determinantal variety
  defined by an~$n\times n$ symmetric  matrix of rank at most~$4$
  whose entries are the traces of simple
  closed curves on the surface.
  This result generalizes the Fricke--Klein--Vogt relations
  to arbitrary $n$, and we determine the structure of the
  kernel of the associated matrix.
\end{abstract}
\subjclass[2000]{}
\maketitle

\section{Introduction}
Let $\mathfrak{X}_{n}$ be the $\mathrm{SL}(2;\mathbb{C})$
character variety of the $n$-punctured sphere
\begin{equation}
  \label{define-X_n}
  \mathfrak{X}_{n}=
  \left.
    \Hom\left(\pi_1(\Sigma_{0,n}), \mathrm{SL}(2;\mathbb{C})\right)
    \middle\slash\!\!\middle\slash
    \mathrm{SL}(2;\mathbb{C})
  \right. ,
\end{equation}
where the fundamental group is
\begin{equation}
  \label{define_pi1}
  \pi_1(\Sigma_{0,n})=
  \left\langle
    C_0, C_1, \dots, C_{n-1}
    ~\middle|~
    C_0 \, C_1 \cdots C_{n-1}=1
  \right\rangle.
\end{equation}
The case  $n=3$ is trivial.
It is well known that
$\mathfrak{X}_4$ is related to the Fricke--Klein--Vogt cubic surfaces
for over a century~\cite{FrickeKleinBook1897}.
In our previous paper~\cite{KHikami26a}, we showed skein-theoretically
that
$\mathfrak{X}_5$ is realized as a determinantal quintic hypersurface.
The purpose of this   note is to
extend this construction to
$\mathfrak{X}_n$ for  $n \geq 5$.

\begin{theorem}
  There is a generically $2$-to-$1$ map
  \begin{equation*}
    \phi: \mathfrak{X}_n \to \mathcal{D}_4 ,
  \end{equation*}
  where the
  determinantal variety~$\mathcal{D}_4$ is defined by
  \begin{equation}
    \label{d4_variety}
    \mathcal{D}_4
    =
    \left\{
      \mathbf{M}
      ~\middle|~
      \rank \mathbf{M} \leq 4
    \right\} ,
  \end{equation}
  and  the symmetric
  matrix~$\mathbf{M}=(M_{a b})_{0 \leq a,b \leq   n-1}$
  is  given by
  \begin{equation*}
    M_{a b}
    =
    \begin{cases}
      2, & \text{for $a=b$},
      \\[4pt]
      \Tr\left(
      \rho(C_a  \, C_{a+1}\cdots C_{b-1})
      \right),
         &
           \text{for $0\leq a<b \leq n-1$},
    \end{cases}
  \end{equation*}
  for $\rho\in \Hom(\pi_1(\Sigma_{0,n}), \mathrm{SL}(2;\mathbb{C}))$.
\end{theorem}
For general properties of determinantal varieties,
see, \emph{e.g.}, \cite{BrunsVette88Book,Dolga12Book}.

\section{Determinantal Variety}
\subsection{Construction of the Gram matrix}

We set
the symmetric  bilinear form  $B(X,Y)$
in the $4$-dimensional vector space $\mathrm{Mat}(2;\mathbb{C})$
by
\begin{equation}
  \label{eq:2}
  \begin{aligned}[b]
    B(X,Y)
    & =
      \det(X+Y)-\det(X)-\det(Y)
    \\
    & =
      \Tr(X) \Tr(Y)-\Tr(X \, Y) .
  \end{aligned}
\end{equation}
Restricting to $\mathrm{SL}(2;\mathbb{C})$,
this gives the non-degenerate bilinear form~$B(X,Y)=\Tr (X^{-1} Y)$.
The isometry group  is $\mathrm{SO}(4;\mathbb{C})$,
and
it  is identified with
the action
\begin{equation*}
\mathrm{SL}(2;\mathbb{C}) \times \mathrm{SL}(2;\mathbb{C})
\cong \mathrm{SO}(4;\mathbb{C}) 
\end{equation*}
on  $n$-tuple $(X_0, X_1, \dots, X_{n-1}) \in \mathrm{SL}(2;\mathbb{C})^n$ by
$X_a \mapsto U \, X_a \, V^{-1}$.
When we set
\begin{equation}
  \rho(C_{a}) =  X_a^{-1} \,  X_{a+1},
\end{equation}
for $0 \leq a \leq n-2$,
the left $\mathrm{SL}(2;\mathbb{C})$ action is free on
$\rho(C_a)$
while the right
action is the conjugaction
$\rho(C_a) \mapsto V \, \rho(C_a) \, V^{-1}$
as in~\eqref{define-X_n}.

The
symmetric Gram matrix~$\mathbf{M}=(M_{ab})_{0\leq a,b <n}$
is then defined by
\begin{equation}
  \label{define_Gram}
  M_{a b}=
  \Tr \left( X_a^{-1} \, X_b \right)
  =
  \begin{cases}
    2, & \text{for $a=b$},
    \\
    \Tr\left(
    \rho(C_a  \, C_{a+1}\cdots C_{b-1})
    \right),
    &
    \text{for $0\leq a<b \leq n-1$}.
  \end{cases}
\end{equation}
We verify that 
$\rank \mathbf{M}=4$ generically,
since
it admits the factorization
\begin{equation}
  \label{eq:1}
  \mathbf{M}
  =
  \mathbf{X}^\top
  \,
  \begin{pmatrix}
    &&& 1\\
    && -1 & \\
    & -1 &&\\
    1 &&&
  \end{pmatrix}
  \,
  \mathbf{X} .
\end{equation}
Here
$\mathbf{X}$ is the $4 \times n$ matrix
$\mathbf{X}=
\left(
  \mathbf{x}_0,\mathbf{x}_1,\dots, \mathbf{x}_{n-1}
\right)$
whose columns are the vectorizations
\begin{equation*}
  \mathbf{x}_a=
  \operatorname{vec}(X_a)=
  \left(x_{11}, x_{21}, x_{12}, x_{22}
  \right)^\top ,
\end{equation*}
where
\begin{equation*}
  X_a=
  \begin{pmatrix}
    x_{11} & x_{12} \\
    x_{21} & x_{22}
  \end{pmatrix}
  \in \mathrm{SL}(2;\mathbb{C}) .
\end{equation*}
By the First  Fundamental Theorem of Invariants for
$\mathrm{SO}(4;\mathbb{C})$~\cite{Weyl39Book},
the invariant ring is generated by the pairings $M_{ab}$ together with
the maximal minors
$\Delta_I=\det(\mathbf{x}_{i_1}, \mathbf{x}_{i_2},
\mathbf{x}_{i_3}, \mathbf{x}_{i_4})$ of $\mathbf{X}$.
The Second Fundamental Theorem provides the  relations among these
generators:
the $5\times 5$ minors of $\mathbf{M}$ vanish,
and
$\Delta_I^2= \det \mathbf{M}_{I I}$.
Consequently,
the map~$\phi$ from the character variety to the determinantal
variety is generically $2$-to-$1$,
in analogy with the case $n=5$~\cite{Simpson18a,KHikami26a}.

The Gram  matrix~$\mathbf{M}$
is a natural generalization of~\cite{KHikami26a}.
For a fixed boundary monodromy
\begin{equation}
  b_a=-\Tr\left(\rho(C_a)\right) ,
\end{equation}
where the minus sign is the convention of  the   skein algebra,
the entries of the  Gram matrix  corresponding to non-boundary simple
closed curves are given  by
\begin{equation}
  p_{a,a+1,\dots, b}
  =-\Tr\left(\rho(C_a \, C_{a+1}\cdots C_b)\right) ,
\end{equation}
where the subscripts
are taken modulo~$n$  with
$2 \leq L=b-a+1\leq
\left\lfloor \frac{n}{2} \right\rfloor  $.
When $n$ is even and $L=\frac{n}{2}$, complementary intervals are
identified, and one representative is chosen from each pair.
This yields
$\frac{1}{2}n(n-3)$ independent coordinates.
The expected  dimension of the symmetric determinantal variety of rank $r=4$
with
fixed~$b_a$ and  fixed diagonal entries is
\begin{equation*}
r \,n-
\begin{pmatrix}
  r \\
  2
\end{pmatrix}
-2n=2n-6,
\end{equation*}
which
agrees with
the expected dimension of~$\mathfrak{X}_n$.

\subsection{Kernel of $\mathbf{M}$}
We recall that
$\mathbf{M}$ is row-equivalent to
\begin{equation*}
  \begin{pmatrix}
    I_4 & \mathbf{M}_4^{-1} \mathbf{m}_4&\cdots
    &
      \mathbf{M}_4^{-1} \mathbf{m}_{n-1}
    \\
    O & O & \cdots & O
  \end{pmatrix} ,
\end{equation*}
where 
we have assumed that the $4 \times 4$ minor
$\mathbf{M}_4=
\left( B(X_a,X_b) \right)_{0\leq a,b \leq 3}$
is non-singular,
and
\begin{equation*}
  \mathbf{m}_a
  =
  \left(
    B(X_0, X_a), B(X_1,X_a), B(X_2,X_a), B(X_3, X_a)
  \right)^\top .
\end{equation*}
To construct the bases of the kernel,
we apply
\begin{equation}
  \label{det=Tr-Tr}
  \det \left(
    \mathbf{x}_a,  \mathbf{x}_b,  \mathbf{x}_c,  \mathbf{x}_d
  \right)
  =
  -\Tr
  \left( X_a^{-1} \, X_b \,  X_c^{-1} \, X_d \right)
  +
  \Tr
  \left( X_a^{-1} \, X_d \, X_c^{-1} \,  X_b \, \right) ,
\end{equation}
where
used is
the vectorization $\mathbf{x}_a=\operatorname{vec}(X_a)$ of
$X_a \in\mathrm{SL}(2;\mathbb{C})$.
One then obtains by the Cramer rule with~\eqref{det=Tr-Tr} that
$\ker \mathbf{M}$ is spanned by the vectors
$\mathbf{v}_k$ for $4\leq k<n$:
\begin{equation}
  \label{eq:3}
  \NiceMatrixOptions{code-for-first-col= {\color{blue}\scriptstyle}}
  \mathbf{v}_k
  =
  \begin{pNiceArray}{cc}[first-col]
    0
    &
    - \Tr\left(
      \rho(C_1 \cdot C_3 C_4\cdots C_{k-1})
    \right)
    + \Tr\left(
      \rho(C_0 C_2\cdot C_{k} C_{k+1}\cdots C_{n-1})
    \right)
    \\[3pt]
    1
    &
    - \Tr\left( \rho(C_2 \cdot C_k C_{k+1}\cdots C_{n-1})
    \right)
    +
    \Tr\left( \rho(C_0 C_1 \cdot C_3 C_4\cdots C_{k-1})
    \right)
    \\[3pt]
    2
    &
    -\Tr\left(  \rho (C_0 \cdot C_3 C_4 \cdots C_{k-1}) \right)
    +
    \Tr\left( \rho(C_1 C_2 \cdot C_k C_{k+1} \cdots C_{n-1})
    \right)
    \\[3pt]
    3&
    - \Tr\left( \rho(C_1\cdot C_k C_{k+1}\cdots C_{n-1})
    \right)
    +
    \Tr\left( \rho(C_0 \cdot C_2 C_3\cdots C_{k-1})
    \right)
    \\
    \Vdots
    &
    0
    \\
    &
    \Vdots
    \\
    &
    0
    \\[3pt]
    k
    &
    - \Tr\left( \rho(C_0 C_2)
    \right)
    + \Tr\left( \rho( C_1 \cdot C_3 C_4 \cdots C_{n-1})
    \right)
    \\
    \Vdots
    &
    0
    \\
    &
    \Vdots
    \\
    n-1
    &
    0
  \end{pNiceArray} .
\end{equation}
In the case  $n=5$, this recovers  the result in~\cite{KHikami26a}.
One may also verify
$\mathbf{v}_k\in \ker \mathbf{M}$ by
skein-theoretic arguments  as in~\cite{KHikami26a}.

\subsection{Action of the Mapping Class Group}

The mapping class group of the $n$-punctured sphere is generated by the
half Dehn twists $\sigma_a$ subject to the
relations~\cite{Birman74,FarbMarg11Book,MagnKarrSoli76Book}
\begin{equation}
  \label{eq:4}
  \Aut(\Sigma_{0,n})
  =
  \left\langle
    \sigma_0, \sigma_1, \dots, \sigma_{n-2}
    \,\middle|\,
    \begin{matrix}
      \sigma_a \, \sigma_{a+1} \, \sigma_a
      =
      \sigma_{a+1} \,  \sigma_a \, \sigma_{a+1}
      \quad \text{for $0\leq a\leq n-3$}
      \\[1pt]
      \sigma_a\, \sigma_b = \sigma_b \, \sigma_a
      \quad \text{for $ |a-b|>1$},
      \\[1pt]
      \sigma_0 \cdots \sigma_{n-3} \, \sigma_{n-2}^2
      \sigma_{n-3} \cdots \sigma_0
      =
      (\sigma_0 \, \sigma_1\cdots \sigma_{n-2})^{n}=1
    \end{matrix}
  \right\rangle .
\end{equation}
The generator $\sigma_a$ acts on the gnerators~\eqref{define_pi1} via 
the Hurwitz action
\begin{equation}
  \label{Hurwitz}
  \sigma_a:
  \left( C_0, \dots, C_a, C_{a+1}, \dots, C_{n-1} \right)
  \mapsto
  \left(
    C_0,  \dots, C_a \, C_{a+1}\,C_a^{-1} , C_a, \dots, C_{n-1}
  \right) .
\end{equation}
One verifies that
$\rank \sigma_a(\mathbf{M})=4$ at a generic point.

\section{Concluding Remarks}

We have identified  the determinantal variety~$\mathcal{D}_4$~\eqref{d4_variety}
as the image of
the character variety of the $n$-punctured sphere under the map
$\phi$.
The singular locus of $\mathcal{D}_4$ consists generically of
$\mathcal{D}_3$
which is the determinantal variety of matrices of rank at most $3$.
It can  have accidental singularities  depending  on the boundary
monodromy $(b_0, b_1, \dots, b_{n-1})$ as in the case
$n=5$~\cite{KHikami26a}.
Future works will
be
to quantize this framework
via  skein algebras, cluster algebras, and the double
affine Hecke algebras, building on our  previous approaches.

\section*{Acknowledgements}
We have used
ChatGPT  and Claude
for communications and 
for   refinement of the exposition of this note.
We take full  responsibility  for the  mathematical content.
The work of KH is supported in part by
JSPS KAKENHI
Grant Number JP23K22388.

\end{document}